\documentclass[a4paper,12pt]{amsart}
\usepackage{amssymb}
\usepackage{ifthen}
\nonstopmode \numberwithin{equation}{section}
\newtheorem{thm}{Theorem}
\newtheorem{cor}{Corollary}
\newtheorem{lem}{Lemma}

\newtheorem{conj}{Conjecture}

\theoremstyle{definition}
\newtheorem{defn}{Definition}[section]

\newtheorem{prob}[equation]{Problem}

\newenvironment{rem}{%
\bigskip
\noindent \textsl{{\sl Remark. }}}{\bigskip}
\newenvironment{rems}{%
\bigskip
\noindent \textsl{{\sl Remarks. }}}{\bigskip}

\newcounter {own}
\def\theown {\thesection       .\arabic{own}}

\newenvironment{pf}[1][]{%
 \vskip 3mm
 \noindent
 \ifthenelse{\equal{#1}{}}%
  {{\slshape Proof. }}%
  {{\slshape #1.} }%
 }%
{\qed\bigskip}

\newcounter{alphabet}
\newcounter{tmp}
\newenvironment{Thm}[1][]{\refstepcounter{alphabet}%
\bigskip%
\noindent%
{\bf Theorem \Alph{alphabet}}%
\ifthenelse{\equal{#1}{}}{}{ (#1)}%
{\bf .} \itshape}{\vskip 8pt}

\newcommand{\ID}{{\mathbb D}}

\newcommand{\IC}{{\mathbb C}}
\newcommand{\C}{{\mathbb C}}
\newcommand{\sphere}{{\widehat{\mathbb C}}}

\newcommand{\inv}{^{-1}}

\newcommand{\dz}{{\partial}}
\newcommand{\dzb}{{\bar\partial}}

\def\be{\begin{equation}}
\def\ee{\end{equation}}

\newcommand{\bee}{\begin{enumerate}}
\newcommand{\eee}{\end{enumerate}}

\newcommand{\blem}{\begin{lem}}
\newcommand{\elem}{\end{lem}}
\newcommand{\bthm}{\begin{thm}}
\newcommand{\ethm}{\end{thm}}
\newcommand{\bcor}{\begin{cor}}
\newcommand{\ecor}{\end{cor}}
\newcommand{\beg}{\begin{examp}}
\newcommand{\eeg}{\end{examp}}
\newcommand{\begs}{\begin{examples}}
\newcommand{\eegs}{\end{examples}}
\newcommand{\bdefe}{\begin{defn}}
\newcommand{\edefe}{\end{defn}}
\newcommand{\bprob}{\begin{prob}}
\newcommand{\eprob}{\end{prob}}
\newcommand{\bei}{\begin{itemize}}
\newcommand{\eei}{\end{itemize}}

\newcommand{\bcon}{\begin{conj}}
\newcommand{\econ}{\end{conj}}
\newcommand{\bcons}{\begin{conjs}}
\newcommand{\econs}{\end{conjs}}
\newcommand{\bprop}{\begin{propo}}
\newcommand{\eprop}{\end{propo}}
\newcommand{\br}{\begin{rem}}
\newcommand{\er}{\end{rem}}
\newcommand{\brs}{\begin{rems}}
\newcommand{\ers}{\end{rems}}
\newcommand{\bo}{\begin{obser}}
\newcommand{\eo}{\end{obser}}
\newcommand{\bos}{\begin{obsers}}
\newcommand{\eos}{\end{obsers}}
\newcommand{\bpf}{\begin{pf}}
\newcommand{\epf}{\end{pf}}
\newcommand{\ba}{\begin{array}}
\newcommand{\ea}{\end{array}}
\newcommand{\beq}{\begin{eqnarray}}
\newcommand{\beqq}{\begin{eqnarray*}}
\newcommand{\eeq}{\end{eqnarray}}
\newcommand{\eeqq}{\end{eqnarray*}}

\newcommand{\ra}{\rightarrow}

\newcounter{minutes}
\divide\time by 60
\newcounter{hours}
\multiply\time by 60 \addtocounter{minutes}{-\time}

\begin{document}
\bibliographystyle{amsplain}
\title{On some results for meromorphic univalent functions having quasiconformal extension}
\begin{center}
{\tiny \texttt{FILE:~\jobname .tex,
        printed: \number\year-\number\month-\number\day,
        \thehours.\ifnum\theminutes<10{0}\fi\theminutes}
}
\end{center}

\author{Bappaditya Bhowmik${}^{~\mathbf{*}}$}
\address{Bappaditya Bhowmik, Department of Mathematics,
Indian Institute of Technology Kharagpur, Kharagpur - 721302, India.}
\email{bappaditya@maths.iitkgp.ernet.in}
\author{Goutam Satpati}
\address{Goutam Satpati, Department of Mathematics,
Indian Institute of Technology Kharagpur, Kharagpur - 721302, India.}
\email{gsatpati@iitkgp.ac.in}

\keywords{ Univalent, Quasiconformal, Distortion estimate, Asymptotic estimate \\
${}^{\mathbf{*}}$ Corresponding author}
\date{ \today ; File: BS.tex}

\begin{abstract}
We consider the class $\Sigma(p)$ of univalent meromorphic
functions $f$ on $\ID$ having simple pole at $z=p\in[0,1)$ with residue 1.
Let $\Sigma_k(p)$ be the class of functions in $\Sigma(p)$ which have $k$-quasiconformal extension to the extended complex plane $\sphere$
where $0\leq k < 1$. We first give a representation formula for functions in this class and using this
formula we derive an asymptotic estimate of the Laurent coefficients for the functions
in the class $\Sigma_k(p)$. Thereafter we give a sufficient condition for functions in  $\Sigma(p)$ to belong in the class $\Sigma_k(p).$
Finally we obtain a sharp distortion result for functions in $\Sigma(p)$ and as a consequence,
we get a distortion estimate for functions in $\Sigma_k(p).$
\end{abstract}

\maketitle
\pagestyle{myheadings}
\markboth{Meromorphic functions having quasiconformal extension}{Meromorphic functions having quasiconformal extension}
\bigskip

\section{Introduction and Preliminary Results}
Let $\IC$ be the complex plane and $\sphere$ be the extended complex plane
$\C\cup\{\infty\}.$ We shall use the following notations throughout the discussion of this article :
$\ID=\{z \in \IC : |z|<1\}$, $\partial\ID=\{z \in \IC : |z|=1\}$, $\overline{\ID}=\{z \in \IC :
|z|\leq 1\}$, $\ID^{*}= \{z \in \IC : |z|>1\}$, $\overline{\ID^{*}}= \{z \in \IC :
|z|\geq1\}$.

The univalent analytic mappings defined in $\ID$ having quasiconformal extension to the whole complex plane play a vital role in Teichm\"{u}ller spaces.
There are number of results for such functions obtained by O. Lehto, R. K\"{u}hnau  and various other mathematicians starting from the work
of L. Ahlfors (see \cite{AB}) in the year 1960 to till date.
We refer to the following articles \cite{Kr1}--\cite{Ku2} for various other results on such mappings.

In this paper our main concern is the univalent meromorphic mappings defined
in $\ID$ with pole at $z=p\in[0,1)$ having quasiconformal extension to the
whole complex plane. O. Lehto extensively studied coefficient problems,
growth estimate for meromorphic functions with pole at the origin $(p=0)$
having quasiconformal extension to the whole plane. We refer to the articles
\cite{Lehto} and \cite{Lehto1} for  further details. In the present
article, we mainly consider the class $\Sigma_k(p)$ of meromorphic univalent
functions with pole at $z=p\in[0,1)$ having quasiconformal extension to the
whole complex plane. This newly defined class of functions has been
introduced and studied in a recent article (compare \cite{BSS}).

Let $f$ be a function with $L^{1}$-derivatives (see \cite[I \S 4.3]{LV1}) in
the whole complex plane $\IC$ and $\Omega$ be a Jordan domain in $\IC$ such
that $\Omega \cup  \partial \Omega=:\overline{\Omega} \subset \IC$ with a
rectifiable boundary curve $\partial \Omega$. We also denote
$\overline{\partial}f:= \partial f/ \partial \overline{z}$ and $\partial f:=
\partial f/\partial z.$ Now applying `Cauchy-Pompeiu' (see \cite[III \S
7]{LV}) formula for such $f$, we get
\begin{equation} \label{eq1}
f(z) = \frac{1}{2\pi i} \int_{\partial \Omega}^{} \frac{f(\zeta)}{\zeta-z} d\zeta - \frac{1}{\pi}\iint\limits_\Omega
\frac{\overline{\partial}f(\zeta)}{\zeta-z} \,d\xi\,d\eta,~~\text{where}~z \in \Omega~~\text{and}~ \zeta=\xi + i \eta.
\end{equation}
If $f(z)\rightarrow 0$ as $z\rightarrow \infty$, then taking
$\Omega=\lbrace\zeta\in\IC : |\zeta|<R\rbrace$, for some $R>0$ and letting $R\to\infty$, the first term of (\ref{eq1}) vanishes and we get,
\begin{equation} \label{eq2}
f(z)= T[\overline{\partial}f](z),
\end{equation}

where
$$
T[\omega](z)=-\frac{1}{\pi} \iint\limits_{\mathbb{C}}
\frac{\omega(\zeta)}{\zeta-z} \,d\xi d\eta.$$ We assume that the function $\omega$ in
above expression belongs to the class $C_{0}^{\infty}$- the class of
infinitely many times differentiable functions with compact support in the
complex plane. We then have

$$
\partial T[\omega](z)=H[\omega](z),
$$
where $H$ is the Hilbert transformation defined by
$$
H[\omega](z):=-\frac{1}{\pi} \iint\limits_{\mathbb{C}}
\frac{\omega(\zeta)}{(\zeta-z)^2} \,d\xi d\eta.
$$
Let $\Sigma$ be the class of univalent meromorphic functions $f$ on $\ID$
having simple pole at the origin with residue 1. Let each $f\in \Sigma$  has the
following expansion
\be\label{eq3a}
f(z)=\frac{1}{z}+\sum_{n=0}^{\infty}b_n z^{n}, \quad z\in \ID.
\ee
In this article our main focus will be the class of function which have pole no more at origin but at a nonzero point. We consider the class $\Sigma(p)$ of univalent meromorphic
functions $f$ on $\ID$ having simple pole at $z=p\in[0,1)$ with residue 1. Therefore, each $f\in \Sigma(p)$  has the following expansion
\be\label{eq3b}
f(z)= \frac{1}{z-p}+ \sum_{n=0}^{\infty}b_n z^{n}, \quad z\in \ID.
\ee
Let  $\Sigma_k$ be the class of functions  in $\Sigma$ that have $k$-quasiconformal extension $(0 \leq k <1)$ to the whole plane $\sphere$ and let
$\Sigma_k(p)$ be the class of functions in $\Sigma(p)$ that have $k$-quasiconformal extension to the whole plane $\sphere$. Here, a mapping $f:\sphere\to\sphere$ is called $k$-quasiconformal
if $f$ is a homeomorphism and has locally $L^2$-derivatives on
$\IC\setminus\{f\inv(\infty)\}$
(in the sense of distribution) satisfying $|\dzb f|\le k|\dz f|$ a.e. .
Thus the complex dilatation $\mu_f(z)$ of $f$ satisfies $|\mu_f(z)|\leq k$, for $z\in\overline{\ID^{*}}$ and vanishes on $\ID$.
Let $\Sigma_k^0(p)$ be the class of functions in  $\Sigma_k(p)$ such that $b_0=0$. Therefore, each $f$ in $\Sigma_k^0(p)$  has the expansion of the
following form

 \begin{equation}\label{eq4}
f(z)= \frac{1}{z-p}+  \sum_{n=1}^{\infty}b_n z^{n}, \quad z\in \ID.
\end{equation}
Let $f\in \Sigma_k^{0}(p)$ and we make a change of variable $\psi(z)=f(1/z)$, so that $\psi$ has the following expansion

 \begin{equation}\label{eq4a}
\psi(z)= z/(1-pz)+  \sum_{n=1}^{\infty}b_n z^{-n}, \quad z\in \ID^{*}.
\end{equation}
As $\psi$ is obtained by composing a M\"{o}bius transformation with a $k$-quasiconformal mapping, therefore it is also a
$k$-quasiconformal map in $\sphere$. Hence the complex dilatation of $\psi$ satisfies $|\mu_{\psi}(z)|\leq k$ for $z \in \overline{\ID}$
and $\mu_{\psi}(z)$ vanishes outside $\overline{\ID}$, i.e. $\mu_{\psi}(z)$ has bounded support. We see that $|\mu_f|=|\mu_{\psi}|=:|\mu|$.
Now, if $f\in\Sigma_k^0(p)$, then $\psi(z)- z/(1-pz) \rightarrow 0$ as $z\rightarrow \infty.$ Thus from (\ref{eq2}), we have

\begin{equation}\label{eq4a1}
 \psi(z)- z/(1-pz) = T\left[ \overline{\partial} \left(\psi(z)- z/(1-pz) \right)\right]= T[\overline{\partial}\psi](z).
\end{equation}
Taking partial derivative of both sides w.r.t. $z$ and using $\partial T [\omega] =H [\omega]$, we get

$$\partial \psi(z) = 1/(1-pz)^2+ H [\overline{\partial}\psi](z).
$$
As $\overline{\partial}\psi= \mu \partial \psi$, the above equation takes the form

\begin{equation}\label{eq5}
\overline{\partial}\psi(z) = \mu/(1-pz)^2 + \mu H [\overline{\partial}\psi](z).
\end{equation}
Now we wish to solve this equation. We will see that it is solvable in $L^2$ but if we assume that $(\|\mu\|_{\infty}\|H\|_q) < 1$, for $q\geq 2$ then it will be solvable in $L^q$.
From `Calderon-Zygmund' inequality (see \cite[I p. 26]{LV1}), we know that $\|H \omega\|_q \leq A_q\|\omega\|_q$, where $A_q$ is a constant, i.e. $\|H\|_q$ is bounded in $L^q.$
In particular we have $\|H\|_2=1.$
We first define inductively
\begin{equation}\label{eq5a}
\phi_1 = \mu/(1-pz)^2,~ \mbox {and}~~~~\phi_n= \mu H [\phi_{n-1}],~~~n=2,3,\cdots .
\end{equation}
Now, since  $\mu(z)=0$ for $z\in\ID^*$, then for $i=2,3,\cdots$, we have
\begin{align}
\|\phi_i\|_q  &= \|\mu H[ \phi_{i-1}]\|_q \nonumber \\ & = |\mu|\|H [\phi_{i-1}]\|_q \nonumber\\
& \leq \|\mu\|_{\infty}\|H\|_q \|\phi_{i-1}\|_q \nonumber\\
&\vdots \nonumber\\
& \leq (\|\mu\|_{\infty}\|H\|_q)^{i-1} \|\phi_1\|_q \label{eq6}.
\end{align}
Next we estimate
\begin{align}
\|\phi_1\|_q &=\bigg( \iint\limits_{|z|\leq 1} \frac{|\mu(z)|^q }{|1-pz|^{2q}} \,dx dy \bigg)^{1/q} \nonumber \\
& \leq \|\mu_{\infty}\| \bigg( \int\limits_{0}^{1} \int\limits_{0}^{2\pi} \frac{r}{(1-pr)^{2q}} \,drd \theta \bigg)^{1/q} \nonumber \\
&= C(p,q)\|\mu\|_{\infty},\label{eq7}
\end{align}
where $C(p,q)$ is a constant depending on $p$ and $q$, and after a little calculation we find it as

\begin{equation}\label{eq7a}
C(p,q)= \left( \frac{2\pi}{p^2}
\left[ \frac{(1-p)^{(2-2q)}}{2-2q} - \frac{(1-p)^{(1-2q)}}{1-2q} + \frac{1}{(1-2q)(2-2q)} \right] \right)^{1/q}.
\end{equation}

\noindent Thus from (\ref{eq6}), we get

\begin{equation}\label{eq8}
\|\phi_i\|_q \leq C(p,q)\|H\|_q^{i-1}\|\mu\|_{\infty}^i, \quad i=2,3,\cdots . 
\end{equation}
Now we define $\omega_n= \sum\limits_{i=1}^{n} \phi_i$ and wish to show that $\lbrace \omega_n \rbrace_{n\geq 1}$ is a Cauchy sequence in $L^q$. Since
$L^{q}$ is complete, $\{\omega_n\}$ will converge in $L^q.$
For $n>m$,
\beqq
\|\omega_n-\omega_m\|_q
&\leq & \sum\limits_{i=m+1}^{n}\|\phi_i\|_q \nonumber\\
& \leq & C(p,q) \sum\limits_{i=m+1}^{n} \|H\|_q^{i-1}\|\mu\|_{\infty}^i \nonumber\\
&=& \frac{C(p,q)}{\|H\|_q} \sum\limits_{i=m+1}^{n} (\|H\|_q\|\mu\|_{\infty})^i \nonumber\\
&\leq & \frac{C(p,q)}{\|H\|_q} \sum\limits_{i=m+1}^{\infty} (\|H\|_q\|\mu\|_{\infty})^i \nonumber\\
&=& \frac{C(p,q)}{\|H\|_q}\left(\frac{M^{m+1}}{1-M}\right) \longrightarrow 0 \quad \mbox{as} \quad n>m \rightarrow \infty,\nonumber
\eeqq
where $M=\|H\|_q||\mu\|_{\infty}<1.$ Hence the sequence $\{\omega_n\}$ is
convergent so that $ \lim\limits_{n\to\infty} \omega_n = \sum\limits_{i=1}^{\infty} \phi_i=: \omega \in L^{q}$.
We now have $H[\omega_n]= \sum\limits_{i=1}^{n} H[\phi_i],$ as $H$ is a linear operator.
Hence,
$$
\mu H [\omega_n]=\sum\limits_{i=1}^{n} \mu H[\phi_i]=\sum\limits_{i=1}^{n} \phi_{i+1}=\sum\limits_{i=2}^{n+1} \phi_i,
$$
which imply
$$
\frac{\mu}{(1-pz)^2} + \mu H [\omega_n] = \sum\limits_{i=1}^{n+1} \phi_i= \omega_{n+1}.
$$
Taking limit both sides of above equation as $n \to \infty$, we have
$$
\frac{\mu}{(1-pz)^2} + \mu H [\omega] = \omega.
$$
So, from above equation it follows that $\omega=\overline{\partial}\psi$ satisfies equation (\ref{eq5}).
 Using this result we provide a representation theorem for functions in $f\in\Sigma_k^0(p)$.
We follow the idea due to Lehto \cite[I \S 4.3]{LV1}.
This is one of the main contents in the next section.

In 1976, J. G. Krzy\.{z} \cite{{Kr76}}, gave a sufficient condition for functions to belong in the class $\Sigma_k(p).$
We state it below.

\begin{Thm}\label{thk}
Let $f \in \Sigma$ have the expansion of the form \eqref{eq3a} in $\ID$.
If there exists $k,~0 \leq k<1$, such that
$$
|z^2f^\prime(z)+1|\leq k|z|^2,\quad{for~all}~z\in\ID.
$$
Then $f\in\Sigma_k.$
\end{Thm}

We also provide a sufficient condition for functions to belong in the class $\Sigma_k(p)$ in the following section.
Next, we state a theorem proved by K. L\"{o}wner \cite{Lowner} for the class $\Sigma$.

\begin{Thm}\label{thl}
Let $f\in\Sigma$ have the expansion of the form \eqref{eq3a}. Then

$$
|z^2f^\prime(z)|\leq\frac{1}{1-|z|^2}, \quad z\in\ID,
$$
where equality holds at a point $z=z_0$ in $\ID$ if and only if
$$
f(z)= \frac{1}{z} + b_0 - \frac{(z_0^{-1}-\overline{z}_0)\overline{z}_0z}{1-\overline{z}_0z}, \quad z\in\ID,
$$
where $b_0$ is a constant.
\end{Thm}

The above theorem  was slightly improved by T. Sugawa \cite[Theorem 1]{Sugawa} as follows:

\begin{Thm}\label{ths}
For $f\in\Sigma$ with the expansion of the form \eqref{eq3a}, the inequality
$$
|z^2f^\prime(z)+ 1|\leq\frac{|z|^2}{1-|z|^2}
$$
holds for each $z\in\ID$. Moreover, equality holds at a point $z=z_0 \in \ID$ if and only if
$$
f(z)= \frac{1}{z} + b_0 - \frac{(z_0^{-1}-\overline{z}_0)\overline{z}_0z}{1-\overline{z}_0z}, \quad z\in\ID,$$
for a constant $b_0$.
\end{Thm}
Using the Area theorem for the class $\Sigma_k$ (see \cite[\S3]{Lehto}) and the above theorem we get that, if $f\in \Sigma_k$ with the expansion of the form \eqref{eq3a} in $\ID$, then
$$|z^2f^\prime(z)+ 1|\leq\frac{k|z|^2}{1-|z|^2}, \quad z\in\ID.
$$
 In the next section we generalize Theorem C for functions in the class $\Sigma(p)$ and as a consequence we obtain a distortion result for functions in $\Sigma_k(p)$.

\section{Main Results}
We start the Section with the following result which we will use to find a representation formula for functions in  $\Sigma_k^0(p)$.

\bthm\label{th1}
If $f\in\Sigma_k^0(p)$, then $f(z)= 1/(z-p)+ \sum\limits_{i=1}^{\infty} T[\phi_i](1/z)$ for $z\in\IC,$ where $\phi_i$'s are defined in $(\ref{eq5a})$.
\ethm

\bpf
Let $\psi(z)=f(1/z), z\in \ID^*$. Therefore from (\ref{eq4a}) we have, $\psi(z)- z/(1-pz) \to 0$ as $z \to \infty$, so by (\ref{eq4a1}) we get
$$
\psi(z) = z/(1-pz) + T [\overline{\partial}\psi](z)=\frac{z}{1-pz}+T\left[\sum\limits_{i=1}^{\infty}\phi_i\right](z).
$$
We have to show that $T\left[\sum\limits_{i=1}^{\infty}\phi_i\right](z)= \sum\limits_{i=1}^{\infty}T[\phi_i](z).$
As $T$ is linear it seems to be obvious but we have to only show that the last series in the above expression is convergent.
Using H\"{o}lder's inequality we now have, ($\zeta=\xi+\mathrm{i}\eta$)
\beqq
|T[\phi_i](z)|&=& \frac{1}{\pi} \left| \iint\limits_{|\zeta| \leq 1} \frac{\phi_i(\zeta)}{\zeta-z}\,d\xi\,d\eta ~\right|\nonumber\\
 & \leq & \frac{1}{\pi} \iint\limits_{|\zeta|\leq 1} \frac{|\phi_i(\zeta)|}{|\zeta-z|}\,d\xi\,d\eta \nonumber\\
&\leq & \frac{1}{\pi} \left( \iint\limits_{|\zeta| \leq 1}| \phi_i(\zeta)|^q \,d\xi\,d\eta \right)^{1/q}\left( \iint\limits_{|\zeta| \leq 1} \frac{1}{|\zeta-z|^s} \,d\xi\,d\eta \right)^{1/s} \quad (1/q+1/s=1)\\
&=& \frac{1}{\pi}\|\phi_i\|_q \Big( \iint\limits_{|\zeta| \leq 1}\frac{1}{|\zeta-z|^s} \,d\xi\,d\eta \Big)^{1/s}.\nonumber
\eeqq
We note that $z$ lies within the whole plane but the integral in the last equation has to be understood by the
Cauchy-Principal value. Now using the estimate derived in (\ref{eq8}), we have
\begin{equation}\label{eq9}
|T[\phi_i](z)| \leq C^{\prime}(p,q)(\|H\|_q\|\mu\|_{\infty})^i,
\end{equation}
where $\|H\|_q$ and $C'(p, q)$ are constants. So applying the Weierstrass-M test, we conclude that the series $\sum\limits_{i=1}^{\infty}T [\phi_i](z)$ is absolutely
and uniformly convergent in $\IC$, and hence we can write

\begin{equation}\label{eq9a}
\psi(z)= \frac{z}{1-pz}+ \sum\limits_{i=1}^{\infty} T[\phi_i](z),\quad z\in\IC.
\end{equation}
Thus we have the following desired representation formula:
$$f(z)= \frac{1}{z-p}+ \sum\limits_{i=1}^{\infty} T[\phi_i](1/z),\quad z\in\IC.
$$
This ends the proof of the Theorem.
\epf

Next, in order to establish asymptotic estimates for $|f(z)-1/(z-p)|$ and $|b_n|$ for functions in the class $\Sigma_k^0(p)$ having expansion (\ref{eq4}),
we first establish another representation formula using previous theorem for functions in the class $\Sigma_k^0(p)$.

\bthm\label{th2}
Let $f\in \Sigma_k^0(p)$ and $k<k_0<1$. As $k\to 0$, we have
$$
f(z)=\frac{1}{z-p}- \frac{1}{\pi} \iint\limits_{\mathbb{D}} \frac{z \mu(\zeta)\,d\xi\,d\eta}{(1-p\zeta)^2(z \zeta-1)} \, + O(k^2),\quad (\zeta=\xi + i \eta)
$$
in the whole plane $\IC$, where $|O(k^2)|\leq ck^2$, and the constant $c$ depends only on $k_0$.
\ethm
\bpf
Since $f\in \Sigma_k^0(p)$ then from (\ref{eq9a}),
we have
$$\psi(z)=z/(1-pz) +T[\phi_1](z)+\sum\limits_{i=2}^{\infty}T[\phi_i](z).
$$
If $q\geq2$ and $k_0\|H\|_q<1$, then $\|\mu\|_{\infty}\|H\|_q\leq k\|H\|_q<k_0\|H\|_q<1.$ Using (\ref{eq9}) we have
\beqq
\sum\limits_{i=2}^{\infty}|T[\phi_i](z)| &\leq & C^{\prime}(p,q)\sum\limits_{i=2}^{\infty}(k\|H\|_q)^i\\
& \leq & C^{\prime}(p,q)k^2\|H\|_q^2 \sum\limits_{i=0}^{\infty}(k_0\|H\|_q)^i \\
&=& k^2\left(\frac{ C^{\prime}(p,q)\|H\|_q^2}{1-k_0\|H\|_q}\right)\\
&=& c k^2,\quad \text{where}~c~\text{depends only on}~k_0.
\eeqq
Hence by the definition of $T$ and $\phi_1$ we have
\begin{equation}\label{eq11}
\psi(z)=\frac{z}{1-pz}- \frac{1}{\pi} \iint\limits_{\mathbb{D}} \frac{\mu(\zeta)d\xi\,d\eta}{(1-p\zeta)^2(\zeta-z)} + O(k^2),\quad  z\in \IC.
\end{equation}
 Hence,
$$
f(z)=\frac{1}{z-p}- \frac{1}{\pi} \iint\limits_{\mathbb{D}} \frac{z \mu(\zeta)d\xi\,d\eta }{(1-p\zeta)^2(z \zeta-1)} \,+ O(k^2),\quad  z\in \IC,
$$
where $|O(k^2)|\leq ck^2$.
\epf

We are now in a position to present the following asymptotic estimate:
\begin{cor}
Each $f\in\Sigma_k^0(p)$ satisfies the following asymptotic bound
\begin{equation}\label{eq12}
\left|f(z)-\frac{1}{z-p}~\right| \leq \frac{k}{\pi} \iint\limits_{\mathbb{D}} \frac{|z|d\xi\,d\eta}{|1-p\zeta|^2|z\zeta-1|} + ck^2,
\,(\zeta=\xi + i \eta), z\in \IC,
\end{equation}
where $|O(k^2)|\leq c k^2$.
\end{cor}
\bpf
It follows from (\ref{eq11}) that
\begin{equation}\label{eq13}
\left|\psi(z)-\frac{z}{1-pz}~\right| \leq \frac{k}{\pi} \iint\limits_{\mathbb{D}} \frac{ d\xi\,d\eta}{|1-p\zeta|^2|\zeta-z|} \,+ ~|O(k^2)|,\quad   z\in \IC.
\end{equation}
Now the inequality (\ref{eq12}) follows by applying a change of variable $f(z)=\psi(1/z)$ in the above inequality. Here we note that in (\ref{eq11}) if
$$
\mu(\zeta)=ke^{i\theta}\frac{\zeta-z}{|\zeta-z|}\Big(\frac{1-p\zeta}{|1-p\zeta|}\Big)^2, \quad\text{a.e. in} ~\ID,
$$
then equality will hold in (\ref{eq13}) and consequently in (\ref{eq12}).
We choose here `$\theta \in (0, 2\pi]$' such that the second and the third term of the right hand side of (\ref{eq11})
have the same argument so that equality holds in (\ref{eq13}).
\epf

\br
We note here that whenever $p\ra 0$ in (\ref{eq13}), we obtain the estimate proved by O. Lehto for functions in $\Sigma_k$ (see f.i. \cite[Cor. 3.2]{LV1}).
\er

Using the above representation formula, we have the following asymptotic coefficient estimates for functions in $\Sigma_k^0(p)$.


\bthm\label{th3}
Let $f(z)\in\Sigma_k^0(p)$, $(0<p<1)$ with the expansion as given in $(\ref{eq4})$. Then

\begin{equation}\label{eq14}
|b_n| \leq 2k\sum_{m=0}^{\infty}\frac{p^{2m}}{n+2m+1} + Ck^2, \quad n\geq 1.
\end{equation}

Here the constant $C$ is given by
$$
C= \frac{C(p)}{(n \pi)^{1/2}(1-k)}, \quad \text{where} \quad C(p)= \left[\frac{\pi(3-p)}{3(1-p)^3}\right]^{1/2}.
$$

Now the equality

$$|b_n| = 2k\sum_{m=0}^{\infty}\frac{p^{2m}}{n+2m+1},
$$

holds for those functions in $\Sigma_k^0(p)$ whose complex dilatation is given by

$$
\mu(z)=k\left(\frac{z}{\overline{z}}\right)^{\frac{n-3}{2}} \left(\frac{z-p}{\overline{z}-p}\right),\quad z\in\ID^*.
$$

\ethm

\bpf
First we note that for $\zeta\in\ID$ and $z\in\ID^*$, we have
$$
\sum\limits_{n=1}^{\infty}\frac{\zeta^{n-1}}{z^n}=-\frac{1}{\zeta-z}.
$$
Therefore, for $z\in\ID^*$, it follows that

\beqq
T[\phi_i](z) &=& -\frac{1}{\pi}\iint\limits_{\mathbb{D}}\frac{\phi_i(\zeta)}{\zeta-z}\,d\xi\,d\eta \\
 &=& \frac{1}{\pi}\iint\limits_{\mathbb{D}}\sum\limits_{n=1}^{\infty}(\phi_i(\zeta) \zeta^{n-1}z^{-n})\,d\xi\,d\eta \\
 &=& \frac{1}{\pi}\sum\limits_{n=1}^{\infty}\Big(\iint\limits_{\mathbb{D}}\phi_i(\zeta) \zeta^{n-1})\,d\xi\,d\eta\Big)z^{-n}.
\eeqq
Now using the representation in (\ref{eq9a}), we get
\beqq
\psi(z)-z/(1-pz)&=& \sum\limits_{i=1}^{\infty}T[\phi_i](z)\\
 &=& \sum\limits_{n=1}^{\infty}\Big(\frac{1}{\pi}\sum\limits_{i=1}^{\infty}\iint\limits_{\mathbb{D}}\phi_i(\zeta) \zeta^{n-1}\,d\xi\,d\eta\Big)z^{-n}.
\eeqq
Thus comparing the above representation of $\psi$ with the expansion in (\ref{eq4a}), we have
the coefficients of $\psi$ as
\begin{equation*}
b_n=\frac{1}{\pi}\sum\limits_{i=1}^{\infty}\iint\limits_{\mathbb{D}}\phi_i(\zeta) \zeta^{n-1}\,d\xi\,d\eta.
\end{equation*}
From the inequality (\ref{eq8}) one can easily obtain
\beqq
\|\phi_i\|_2 & \leq & C(p,2)\,\|H\|_2^{i-1}\|\mu\|_\infty^i \\
& \leq & C(p)k^i,
\eeqq
by virtue of the fact that $\|H\|_2=1$ and $\|\mu\|_\infty \leq k$ and denoting $C(p,2)=:C(p).$ Here we get the value of $C(p)$ by putting $q=2$ in (\ref{eq7a}) as
$$
C(p)=  \left[\frac{\pi(3-p)}{3(1-p)^3}\right]^{1/2}.
$$
Now applying Cauchy-Schwartz inequality in $L^2$, we obtain for each $i\geq 2$,

\begin{align}
\left|\iint\limits_{\mathbb{D}}\phi_i(\zeta) \zeta^{n-1}\,d\xi\,d\eta \right|
&\leq  \iint\limits_{\mathbb{D}}|\phi_i(\zeta)| ~|\zeta|^{n-1}\,d\xi\,d\eta \nonumber\\
&\leq  \|\phi_i\|_2~ \left(\iint\limits_{\mathbb{D}} |\zeta|^{2(n-1)}\,d\xi\,d\eta  \right)^{1/2} \nonumber\\
&=  \|\phi_i\|_2~(\pi/n)^{1/2}\nonumber\\
&\leq  C(p)k^i(\pi/n)^{1/2}.\label{eq14a}
\end{align}
\newline
Now we can write

\begin{equation*}
b_n=\frac{1}{\pi}\iint\limits_{\mathbb{D}}\frac{\mu(\zeta)\zeta^{n-1}}{(1-p\zeta)^2}\,d\xi\,d\eta+\frac{1}{\pi}\sum\limits_{i=2}^{\infty}\iint\limits_{\mathbb{D}}\phi_i(\zeta) \zeta^{n-1}\,d\xi\,d\eta.
\end{equation*}
Hence using (\ref{eq14a}) we get,

\begin{align}
\left|\frac{1}{\pi}\sum\limits_{i=2}^{\infty}\iint\limits_{\mathbb{D}}\phi_i(\zeta) \zeta^{n-1}\,d\xi\,d\eta \right|  &\leq \frac{C(p)}{\pi} \left(\frac{\pi}{n}\right)^{1/2} \sum\limits_{i=2}^{\infty}k^i \nonumber \\
&= \frac{C(p)}{(n\pi)^{1/2}}\left(\frac{k^2}{1-k}\right) \nonumber \\
&= Ck^2,\quad \text{where}\quad C = \frac{C(p)}{(n \pi)^{1/2}(1-k)}.\label{eq14b}
\end{align}
Consequently, we have the following asymptotic representation of the coefficients of the functions in the class $\Sigma_k^0(p)$:

\begin{equation*}
b_n=\frac{1}{\pi}\iint\limits_{\mathbb{D}}\frac{\mu(\zeta)\zeta^{n-1}}{(1-p\zeta)^2}\,d\xi\,d\eta+O(k^2), \quad n=1,2,\cdots.
\end{equation*}
where $|O(k^2)|\leq Ck^2$. Hence we have,
\beqq
|b_n|&\leq & \frac{k}{\pi}\iint\limits_{\mathbb{D}}\frac{|\zeta|^{n-1}}{|1-p\zeta|^2}\,d\xi\,d\eta+Ck^2\nonumber\\
&=& \frac{k}{\pi}\int_{0}^{2\pi}\int_{0}^{1}\frac{r^{n}}{1-2pr\cos\theta+p^2r^2}\,dr\,d\theta~+Ck^2. \nonumber \\
\eeqq
To find the value of the last integral we use the fact that

$$
\frac{1}{2\pi}\int\limits_{0}^{2\pi}\frac{1-r^2}{1-2r\cos\theta +r^2}\,d\theta =1,\quad \text{for}~r<1.
$$
Using this formula and noting that $pr<1$, we get from above

\beqq
|b_n| &\leq & \frac{k}{\pi}\int\limits_{0}^{1} \frac{r^n}{1-p^2r^2} \left(\int\limits_{0}^{2\pi}\frac{1-p^2r^2}{1-2pr\cos \theta +p^2r^2}\,d\theta \right)\,dr+Ck^2\\
&=& 2k\int\limits_{0}^{1}\frac{r^n}{1-p^2r^2}\,dr+Ck^2\\
&=& 2k\int\limits_{0}^{1}r^n(1+p^2r^2+p^4r^4+p^6r^6+\cdots)\,dr+Ck^2\\
&=& 2k\sum_{m=0}^{\infty}\frac{p^{2m}}{n+2m+1}+Ck^2.
\eeqq

%
%
%
\noindent Now the following equality
$$
|b_n|= 2k\sum_{m=0}^{\infty}\frac{p^{2m}}{n+2m+1}
$$
holds whenever
$$
|b_n|=\frac{k}{\pi}\iint\limits_{\mathbb{D}}\frac{|\zeta|^{n-1}}{|1-p\zeta|^2}\,d\xi\,d\eta .
$$
Therefore, it is clear that the above equality holds for the functions whose complex dilatation is given by

$$
\mu(\zeta)=k\left(\overline{\zeta}/\zeta\right)^{\frac{n-1}{2}}\left(\frac{(1-p\zeta)}{|1-p\zeta|}\right)^2 \quad\text{a.e.~in}~ \ID,
$$

i.e.,
$$
\mu(z)=k\left(z / \overline{z}\right)^{\frac{n-3}{2}} \left(\frac{z-p}{\overline{z}-p}\right)\quad\text{for}~z\in\ID^*.
$$

\epf

%

\br
For the case $p=0$ i.e. if $f\in \Sigma_k$ with the expansion as given by (\ref{eq3a}) in $\ID$, the estimate given by (\ref{eq14}) becomes
$$
|b_n| \leq \frac{2k}{n+1} + Ck^2.
$$
In this case the constant in (\ref{eq7}) is given by $C(q)= \pi^{1/q}$ so that the constant in (\ref{eq14a}) is replaced by $C(p)=\pi^{1/2}$.
Hence the constant in (\ref{eq14b}) finally coming as $C=\frac{n^{-1/2}}{1-k}$. This result was proved by O. Lehto (see \cite[II p.74]{LV1}).
\er

Next, we prove  a sufficient condition for functions to belong in the class $\Sigma_k(p).$

\bthm
Let $f \in \Sigma(p)$ has an expansion of the form \text{$(\ref{eq3b})$} in $\ID$. If there exists $k,0 \leq k<1$ and $p,~0\leq p<1$, such that
$$
|(z-p)^2f^\prime(z)+ 1 |\leq \frac{k|z-p|^2}{(1+p)^2} \quad\text{for~all}~z\in\ID,
$$
then $f \in \Sigma_k(p)$.
\ethm
\bpf
For $z\in\ID$, we have
$$
f(z)= 1/(z-p) +  \sum_{n=0}^{\infty}b_n z^{n}.
$$
Therefore we can write
$$
f(z)= 1/(z-p) + \omega(z),
$$
where $\omega(z)=\sum\limits_{n=0}^{\infty}b_nz^n$ is analytic in $\ID$. Hence for $z\in\ID$,
$$
f^\prime(z)+ (z-p)^{-2}= \omega^\prime(z).
$$
It follows from the given condition that $|\omega^\prime(z)|\leq \frac{k}{(1+p)^2}$, $z\in\ID$.
Hence by \\ \cite[Theorem 2]{BSS}, we conclude that $f \in\Sigma_k(p)$.
\epf

%
%

Our next result deals with a distortion inequality for functions in the class $\Sigma(p)$. 

\bthm\label{th5}

Each function $f\in\Sigma(p)$ of the form \text{$(\ref{eq3b})$} satisfies the inequality
\begin{equation}\label{eq16}
\left|f^{\prime}(z)+\frac{1}{(z-p)^2}\right| \leq \frac{1}{(1-p^2)(1-|z|^2)},\quad\text{for}\quad z\in\ID.
\end{equation}
Moreover, equality holds for a constant $b_0$ at $z=z_0\in\ID$, if and only if
\begin{equation}\label{eq17}
f(z)=\frac{1}{z-p} + b_0 - \left(\frac{z_0^{-1}-\overline{z}_0}{1-p^2}\right) \left( \frac{\overline{z}_0-p}{z_0-p}\right) \left(\frac{\overline{z}_0z}{1-\overline{z}_0z}\right).
\end{equation}
\ethm

\bpf
Let $f\in\Sigma(p)$ with the expansion as
$$
f(z)= \frac{1}{z-p}+  \sum_{n=0}^{\infty}b_n z^{n},\quad\text{for} \quad z\in \ID.
$$
We now give a change of variable $z=1/ \zeta$ and let $f(z)=f(1/\zeta)=\psi(\zeta).$ So $\psi$ is defined in $\ID^{*}$ and has Laurent's series expansion as
\beq\label{eq18}
\psi(\zeta)= \frac{\zeta}{1-p \zeta}+  \sum_{n=0}^{\infty}b_n \zeta^{-n}, \quad \zeta\in \ID^{*}.
\eeq
Thus by Chichra's area theorem (see \cite{Chichra}), we have
$$
\sum_{n=1}^{\infty}n|b_n|^2\leq \frac{1}{(1-p^2)^2}.
$$
From equation (\ref{eq18}), taking derivative of both sides we get,
$$
\psi^{\prime}(\zeta)- \frac{1}{(1-p \zeta)^2}  = \sum_{n=1}^{\infty}(-n b_n) \zeta^{-n-1},\quad \zeta\in \ID^{*}.
$$
Now applying Cauchy-Schwartz inequality in above we have,

\begin{align}\label{eq18a}
\left|\psi^{\prime}(\zeta)-\frac{1}{(1-p\zeta)^2}\right| & \leq   \sum_{n=1}^{\infty} n|b_n||\zeta|^{-n-1} \nonumber\\
 & \leq  \sqrt{\sum_{n=1}^{\infty}n|b_n|^2 \sum_{n=1}^{\infty}n|\zeta|^{-2n-2}} \nonumber\\
 & =  \sqrt{\sum_{n=1}^{\infty}n|b_n|^2\frac{1}{(|\zeta|^2-1)^2}}.
\end{align}

Now using Chichra's area theorem in above inequality, we get
\beq\label{eq19}
\left|\psi^{\prime}(\zeta)-\frac{1}{(1-p\zeta)^2}\right| \leq \frac{1}{(1-p^2)(|\zeta|^2-1)}, \quad\text{for}~ \zeta\in \ID^{*}.
\eeq
Returning back to the original variable $z=1/\zeta$ and noting that $-z^2f^\prime(z)=\psi^\prime(\zeta)$, the above inequality yields
$$
\left|z^2f^{\prime}(z)+\frac{z^2}{(z-p)^2}\right| \leq \frac{|z|^2}{(1-p^2)(1-|z|^2)},\quad\text{for}\quad z\in\ID,
$$
which gives us the desired result.

Next we consider the equality case. Now equality holds for $f(z)$ in (\ref{eq16}) for some point $z=z_0\in\ID$ if and only if it does hold  for $\psi(\zeta)$ in (\ref{eq19}) for the point $\zeta=\zeta_0=1/z_0\in \ID^{*}$. This implies that equality will occur in Cauchy-Schwartz inequality as well as in the Chichra's area theorem. For the first case we can say that the two sequences $\lbrace b_n \rbrace$ and $\lbrace (\overline{\zeta}_0)^{-n-1}\rbrace$ are proportional, which means there exists a complex constant $a$ such that $b_n=a(\overline{\zeta}_0)^{-n}$ for $n\geq1$. Again, if equality holds in the Chichra's area theorem, then we have $|a|=\frac{(R-1/R)}{1-p^2}$, where $R=|\zeta_0|>1$. Hence

\begin{eqnarray}\label{eq20}
\psi(\zeta)& = & \frac{\zeta}{1-p\zeta} + b_0 + a\sum_{n=1}^{\infty}(\overline{\zeta}_0)^{-n}\zeta^{-n} \nonumber\\ &=& \frac{\zeta}{1-p\zeta}+b_0+\frac{a}{\overline{\zeta}_0\zeta-1}, \quad\text{for}~ \zeta\in \ID^{*}.
\end{eqnarray}
For this function clearly equality holds in (\ref{eq19}) at the point $\zeta=\zeta_0\in \ID^*$ and hence it does hold in (\ref{eq16}) at $z_0=1/\zeta_0\in \ID$. Next, we show that $\psi$ in (\ref{eq20}) is univalent in $\ID^*$ for a particular value of $a$.

We first make a change of variable $\eta= \frac{\zeta-p}{1-p\zeta}$. Then $\eta\in\ID^*$ if and only if $\zeta\in\ID^*$. Using this transformation and letting $b=\overline{\zeta}_0$, we get from (\ref{eq20})
\beqq
(1-p^2)\psi\left(\frac{\eta+p}{1+p\eta}\right) & = & \eta + \frac{a(1-p^2)}{b(\frac{\eta+p}{1+p\eta})-1} + K,\quad\text{where}~K~ \text{is a constant.}\\
&=& \eta + \frac{a\left(\frac{1-p^2}{1-bp}\right)}{(\frac{b-p}{1-bp})\eta-1}+ \frac{ap\left(\frac{1-p^2}{1-bp} \right)\eta}{(\frac{b-p}{1-bp})\eta-1} - \frac{ap(1-p^2)}{b-p}\\
&& \quad + ~K_0,~ K_0~\text{is another constant.} \\
&=& \eta +  \frac{a\left(\frac{1-p^2}{1-bp}\right)}{(\frac{b-p}{1-bp})\eta-1} + \frac{ap(1-p^2)}{(b-p)\left((\frac{b-p}{1-bp})\eta-1\right)}+ K_0\\
& = & \eta + \frac{\frac{ab(1-p^2)^2}{(1-bp)(b-p)}}{( \frac{b-p}{1-bp})\eta-1}+K_0 \\
& = & \eta + \frac{A}{B\eta -1} + K_0~, \\
\eeqq
where
$$
A=\frac{ab(1-p^2)^2}{(1-bp)(b-p)}\quad\text{and}\quad B=\frac{b-p}{1-bp}.
$$
Let us denote
$$
\phi(\eta)=(1-p^2)\psi\left(\frac{\eta+p}{1+p\eta}\right),\quad~ \eta\in \ID^{*},
$$
then from the above calculation we see that $\phi$ takes the form

\begin{equation}\label{eq21}
\phi(\eta)= \eta + \frac{A}{B\eta -1} + K_0~, \quad~ \eta\in \ID^{*}.
\end{equation}
Hence it is clear that $\phi$ is univalent in $\ID^*$ if and only if $\psi$
is univalent in $\ID^*$. Now we apply  \cite[Lemma 2]{Sugawa} in (\ref{eq21})
to obtain that $\phi$ is univalent in $\ID^*$ precisely when
$|A+\overline{B}-B^{-1}|+|AB| \leq |B^2|-1. $ Using the relation
$|a|=(|b|-|b|^{-1})/(1-p^2)$ and performing a simple calculation, the above
inequality gives
$$
a= -\left(\frac{\zeta_0-(\overline{\zeta}_0)^{-1}}{1-p^2}\right) \left(\frac{1-p\overline{\zeta}_0}{1-p\zeta_0}\right).
$$
Putting this value of $a$ in (\ref{eq20}), we get the univalent extremal function for the inequality (\ref{eq19}) at the point $\zeta=\zeta_0 \in \ID^*$ as
$$
\psi(\zeta)= \frac{\zeta}{1-p\zeta} + b_0 -\left(\frac{\zeta_0-(\overline{\zeta}_0)^{-1}}{1-p^2}\right) \left(\frac{1-p\overline{\zeta}_0}{1-p\zeta_0}\right) \left(\frac{1}{\overline{\zeta}_0\zeta-1} \right), \quad~\zeta\in \ID^*.
$$
Consequently, a change of variable $f(z)=\psi(1/\zeta)$ will yield the required extremal function (\ref{eq17}) for the inequality (\ref{eq16}).
\epf


%
Using the area theorem for $\Sigma_k(p)$ (see \cite[Theorem 1]{BSS}) in $(\ref{eq18a})$ and the above theorem,
we now have the following distortion result. However, sharpness of this bound is not being established.

\begin{cor}
Let  $f\in \Sigma_k(p)$ and have the expansion \eqref{eq3b} in $\ID$. Then
\begin{equation*}
\left|f^{\prime}(z)+\frac{1}{(z-p)^2}\right| \leq \frac{k}{(1-p^2)(1-|z|^2)},\quad\text{for}\quad z\in\ID.
\end{equation*}
\end{cor}
\vspace{1cm}
\noindent {\bf Acknowledgement:}
The authors thank Toshiyuki Sugawa for his suggestions and careful reading of the manuscript.


\begin{thebibliography}{}
\bibitem{AB} {\sc L. Ahlfors} and {L. Bers}, \textit {Riemann's mapping theorem for variable metrics}, Ann. of Math. (2) 72, 385--404, 1960.
\bibitem{LV} {\sc O. Lehto} and {\sc K. I. Virtanen}, \textit{Quasiconformal Mappings in the Plane, 2nd Ed.,} Springer-Verlag, 1973.
\bibitem{LV1} {\sc O. Lehto}, \textit{Univalent Functions and Teichmuller Spaces,} Springer-Verlag, 1987.
\bibitem{Chichra} {\sc P. N. Chichra,} \textit{An area theorem for bounded univalent functions,} Proc. Camb. Phil. Soc., 66, 317-321, 1969.
\bibitem{Lehto} {\sc O. Lehto,} \textit{Schlicht functions with a quasiconformal extension,} Ann. Acad. Sci. Fenn. Ser. A.I., 500, 3-10, 1971.
\bibitem{Lehto1} {\sc O. Lehto,} \textit{Quasiconformal mappings and singular integrals,} Instit. Naz. Alta Mat., Simposia Mathematica, Academic Press, London, 18, 429-453, 1976.
\bibitem{Kr76} {\sc J. G. Krzy\.{z},} \textit{Convolution and quasiconformal extension,} Comment. Math. Helv., 51, 99-104, 1976.
\bibitem{BSS} {\sc B. Bhowmik}, {\sc G. Satpati} and {\sc T. Sugawa,} \textit{Quasiconfornal extension of meromorphic functions with nonzero pole,} Proc. Amer. Math. Soc., 144, 2593-2601, 2016.
\bibitem{Lowner} {\sc K. L\"{o}wner,} \textit{\"{U}ber Extremums\"{a}tze bei der konformen Abbildung des \"{A}u{\ss}eren des Einheitskreises,} Math. Z., 3, 65-67, 1919.
\bibitem{Sugawa} {\sc T. Sugawa,} \textit{A Remark on L\"{o}wner's Theorem,} Inter. Inf. Sc., 18(1), 19-22, 2012.
\bibitem{Kr1} {\sc S. Krushkal,} \textit{General distortion theorem for univalent functions with quasiconformal extension,}
Complex Anal. Oper. Theory, DOI:10.1007/s11785-016-0614-8.
\bibitem{Kr2} {\sc S. Krushkal,} \textit{On the problem of the coefficients for univalent functions with a quasiconformal continuation,} Dokl. Akad. Nauk. SSSR, 287, No. 3, 547-550, 1986.
\bibitem{Kr3} {\sc S. Krushkal,} \textit{Exact coefficient estimates for univalent functions with quasiconformal extension,} Ann. Acad. Sci. Fenn. Ser. A.I. Math., 20, 349–357, 1995.
\bibitem{Kr4} {\sc S. Krushkal,} \textit{Some extremal problems for schlicht analytic functions,} Soviet Math. Dokl., 9, 1968.
\bibitem{Ku1} {\sc R. K\"{u}hnau,} \textit{Wertannahmeprobleme bei quasikonformen Abbildungen mit ortsabh\"{a}ngiger Dilatationsbeschr\"{a}nkung,} Math. Nachr., 40, 1-11, 1969.
\bibitem{Ku2} {\sc R. K\"{u}hnau,} \textit{Verzerrungss\"{a}tze und Koeffizientenbedingungen vom Grunskyschen Typ f\"{u}r quasikonforme Abbildungen,} Math. Nachr., 48, 77-105, 1971.
\end{thebibliography}
\end{document}